\documentclass[10pt,a4paper]{amsart}

\theoremstyle{plain}
\newtheorem{theorem}{Theorem}[section]

\newtheorem{corollary}[theorem]{Corollary}
\theoremstyle{definition}

\theoremstyle{remark}
\newtheorem{remark}[theorem]{Remark}

\def\bin #1#2 {\left( \matrix { #1 \cr #2 \cr } \right) }

\begin{document}

\title[A lower bound for $\chi (\mathcal O_S)$]
{A lower bound for $\chi (\mathcal O_S)$}

\author{Vincenzo Di Gennaro }
\address{Universit\`a di Roma \lq\lq Tor Vergata\rq\rq, Dipartimento di Matematica,
Via della Ricerca Scientifica, 00133 Roma, Italy.}
\email{digennar@axp.mat.uniroma2.it}

\abstract Let $(S,\mathcal L)$ be a smooth, irreducible, projective,
complex surface, polarized by a very ample line bundle $\mathcal L$
of degree $d > 25$. In this paper we  prove that $\chi (\mathcal
O_S)\geq -\frac{1}{8}d(d-6)$. The bound is sharp, and $\chi
(\mathcal O_S)=-\frac{1}{8}d(d-6)$ if and only if $d$ is even, the
linear system $|H^0(S,\mathcal L)|$ embeds $S$ in a smooth rational
normal scroll $T\subset \mathbb P^5$ of dimension $3$, and here, as
a divisor, $S$ is linearly equivalent to $\frac{d}{2}Q$, where $Q$
is a quadric on $T$. Moreover, this is equivalent to the fact that
the general hyperplane section $H\in |H^0(S,\mathcal L)|$ of $S$ is
the projection of a curve $C$ contained in the Veronese surface
$V\subseteq \mathbb P^5$, from a point $x\in V\backslash C$.



\endabstract
\maketitle

\bigskip
\section{Introduction}

In \cite{DGF}, one proves a sharp lower bound for the
self-intersection $K^2_S$ of the canonical bundle of a smooth,
projective, complex surface $S$, polarized by a very ample line
bundle $\mathcal L$, in terms of its degree
$d={\text{deg}}\,\mathcal L$, assuming $d>35$. Refining the line of
the proof in \cite{DGF}, in the present paper we deduce a similar
result for the Euler characteristic $\chi(\mathcal O_S)$ of $S$
\cite[p. 2]{BV}, in the range $d>25$. More precisely, we prove the
following:

\begin{theorem}\label{lbound}  Let $(S,\mathcal L)$ be a smooth,
irreducible, projective, complex surface, polarized by a very ample
line bundle $\mathcal L$ of degree $d > 25$. Then:
$$
\chi (\mathcal O_S)\geq -\frac{1}{8}d(d-6).
$$
The bound is sharp, and the following properties are equivalent.

\medskip
(i) $\chi (\mathcal O_S)= -\frac{1}{8}d(d-6)$;

\medskip
(ii) $h^0(S,\mathcal L)=6$, and the linear system $|H^0(S,\mathcal
L)|$ embeds $S$ in $\mathbb P^5$ as a scroll with sectional genus
$g=\frac{1}{8}d(d-6)+1$;

\medskip
(iii) $h^0(S,\mathcal L)=6$, $d$ is even, and the linear system
$|H^0(S,\mathcal L)|$ embeds $S$  in a smooth rational normal scroll
$T\subset \mathbb P^5$ of dimension $3$, and here $S$ is linearly
equivalent to $\frac{d}{2}(H_T-W_T)$, where $H_T$ is the hyperplane
class of $T$, and $W_T$ the ruling (i.e. $S$ is linearly equivalent
to an integer multiple of a smooth quadric $Q\subset T$).
\end{theorem}

By Enriques' classification, one knows that if $S$ is unruled or
rational, then $\chi (\mathcal O_S)\geq 0$. Hence, Theorem
\ref{lbound} essentially  concerns  irrational ruled surfaces.

In the range $d>35$, the family of extremal surfaces for $\chi
(\mathcal O_S)$ is exactly the same for $K^2_S$. We point out there
is a relationship between this family and the Veronese surface. In
fact one has the following:

\begin{corollary}\label{Veronese}
Let $S\subseteq \mathbb P^r$ be a nondegenerate, smooth,
irreducible, projective, complex surface, of degree $d > 25$. Let
$L\subseteq \mathbb P^r$ be a general hyperplane. Then
$\chi(\mathcal O_S)=-\frac{1}{8}d(d-6)$ if and only if $r=5$, and
there is a curve $C$ in the Veronese surface $V\subseteq \mathbb
P^5$ and a point $x\in V\backslash C$ such that the general
hyperplane section $S\cap L$ of $S$ is the projection
$p_x(C)\subseteq L$ of $C$ in $L\cong\mathbb P^4$, from the point
$x$.
\end{corollary}

In particular, $S\cap L$ is not linearly normal, instead $S$ is.

\bigskip
\section{Proof of Theorem \ref{lbound}}

\begin{remark}\label{k8} $(i)$  We say that $S\subset \mathbb P^r$ is a {\it
scroll} if $S$ is a  $\mathbb P^1$-bundle over a smooth curve, and
the restriction of $\mathcal O_S(1)$ to a fibre is $\mathcal
O_{\mathbb P^1}(1)$. In particular, $S$ is a geometrically ruled
surface, and therefore $\chi (\mathcal O_S)= \frac{1}{8}K^2_S$
\cite[Proposition III.21]{BV}.

\smallskip $(ii)$ By Enriques' classification \cite[Theorem X.4 and Proposition
III.21]{BV}, one knows that if $S$ is unruled or rational, then
$\chi (\mathcal O_S)\geq 0$, and if $S$ is ruled with irregularity
$>0$, then  $\chi (\mathcal O_S)\geq \frac{1}{8}K^2_S$. Therefore,
taking into account previous remark, when $d>35$, Theorem
\ref{lbound} follows from \cite[Theorem 1.1]{DGF}. In order to
examine the range $25<d\leq 35$, we are going to refine the line of
the argument in the proof of \cite[Theorem 1.1]{DGF}.

\smallskip $(iii)$ When $d=2\delta$ is even, then $\frac{1}{8}d(d-6)+1$ is the genus of
a plane curve of degree $\delta$, and the genus of a curve of degree $d$ lying on
the Veronese surface.
\end{remark}

Put $r+1:=h^0(S,\mathcal L)$. Therefore,  $|H^0(S,\mathcal L)|$
embeds $S$ in $\mathbb P^r$. Let $H\subseteq \mathbb P^{r-1}$ be the
general hyperplane section of $S$, so that $\mathcal L\cong \mathcal
O_S(H)$. We denote by $g$ the genus of $H$. If $2\leq r\leq 3$, then
$\chi (\mathcal O_S)\geq 1$. Therefore, we may assume $r\geq 4$.

\bigskip
{\bf The case $r=4$}.

\smallskip
We first examine the case $r=4$. In this case we only have to prove
that, for $d>25$, one has $\chi (\mathcal O_S)> -\frac{1}{8}d(d-6)$.
We may  assume that $S$ is an irrational ruled surface, so
$K^2_S\leq 8\chi (\mathcal O_S)$ (compare with previous Remark
\ref{k8}, $(ii)$). We argue by contradiction, and assume also that
\begin{equation}\label{sest}
\chi (\mathcal O_S)\leq -\frac{1}{8}d(d-6).
\end{equation}
We are going to prove that this assumption implies $d\leq 25$, in contrast with our
hypothesis $d>25$.

By the double point formula:
$$
d(d-5)-10(g-1)+12\chi(\mathcal O_S)=2K^2_S,
$$
and $K^2_S\leq 8\chi (\mathcal O_S)$, we get:
$$
d(d-5)-10(g-1)\leq 4\chi(\mathcal O_S).
$$
And from $\chi (\mathcal O_S)\leq -\frac{1}{8}d(d-6)$ we obtain
\begin{equation}\label{fest}
10g\geq \frac{3}{2}d^2-8d+10.
\end{equation}

\smallskip

Now we distinguish two cases, according that $S$ is not
contained in a hypersurface of degree $<5$ or not.

First suppose that $S$ is not contained in a hypersurface of
$\mathbb P^4$ of degree $<5$. Since $d > 16$, by Roth's Theorem
(\cite[p. 152]{R}, \cite[p. 2, (C)]{EP}),  $H$ is not contained in a
surface of $\mathbb P^3$ of degree $<5$. Using Halphen's bound
\cite{GP}, we deduce that
$$
g\leq \frac{d^2}{10} +
\frac{d}{2}+1-\frac{2}{5}(\epsilon+1)(4-\epsilon),
$$
where $d-1=5m+\epsilon$, $0\leq \epsilon<5$. It follows that
$$
\frac{3}{2}d^2-8d+10\leq \,10 g\,\leq
d^2+5d+10\left(1-\frac{2}{5}(\epsilon+1)(4-\epsilon)\right).
$$
This implies that $d\leq 25$, in contrast with our hypothesis $d>25$.

\smallskip
In the second case,  assume that $S$ is contained in an irreducible
and reduced hypersurface of degree $s\leq 4$. When $s\in\{2,3\}$,
one knows that, for $d>12$,  $S$ is of general type \cite[p.
213]{BF}. Therefore, we only have to examine the case $s=4$. In this
case $H$ is contained in a surface of $\mathbb P^3$ of degree $4$.
Since $d>12$, by Bezout's Theorem, $H$ is not contained in a surface
of $\mathbb P^3$ of degree $<4$. Using Halphen's bound \cite{GP},
and \cite[Lemme 1]{EP}, we get:
$$
\frac{d^2}{8}-\frac{9d}{8}+1\leq \, g\,\leq \frac{d^2}{8}+1.
$$
Hence, there exists a rational number $0\leq x\leq 9$   such that
$$
g=\frac{d^2}{8}+d\left(\frac{x-9}{8}\right)+1.
$$

If $0\leq x\leq \frac{15}{2}$, then $g\leq
\frac{d^2}{8}-\frac{3}{16}d+1$, and from (\ref{fest}) we get
$$
\frac{3}{20}d^2-\frac{4}{5}d+1\,\leq g\,\leq
\frac{d^2}{8}-\frac{3}{16}d+1.
$$
It follows $d\leq 24$,  in contrast with our hypothesis $d>25$.

Assume $\frac{15}{2}< x\leq 9$. Hence,
$$
\left(\frac{d^2}{8}+1\right)-g=
-d\left(\frac{x-9}{8}\right)<\frac{3}{16}d.
$$
By \cite[proof of Proposition 2, and formula (2.2)]{D}, we have
$$
\chi(\mathcal O_S)\geq 1+
\frac{d^3}{96}-\frac{d^2}{16}-\frac{5d}{3}-\frac{349}{16}-(d-3)\left[\left(\frac{d^2}{8}+1\right)-g\right]
$$
$$
>1+
\frac{d^3}{96}-\frac{d^2}{16}-\frac{5d}{3}-\frac{349}{16}-(d-3)\frac{3}{16}d
= \frac{d^3}{96}-\frac{d^2}{4}-\frac{53}{48}d-\frac{333}{16}.
$$
Combining with (\ref{sest}), we get
$$
\frac{d^3}{96}-\frac{d^2}{4}-\frac{53}{48}d-\frac{333}{16}+\frac{1}{8}d(d-6)<0,
$$
i.e.
$$
d^3-12d^2-178d-1998<0.
$$
It follows $d\leq 23$,  in contrast with our hypothesis $d>25$.

This concludes the analysis of the case $r=4$.

\bigskip
{\bf The case $r\geq 5$}.

\smallskip
When $r\geq 5$, by \cite[Remark 2.1]{DGF},  we know that, for $d>5$,
one has $K^2_S>-d(d-6)$, except when $r=5$, and the surface $S$ is a
scroll, $K^2_S=8\chi (\mathcal O_S)=8(1-g)$, and
\begin{equation}\label{bound}
g=\frac{1}{8}d^2-\frac{3}{4}d+\frac{(5-\epsilon)(\epsilon+1)}{8},
\end{equation}
with $d-1=4m+\epsilon$, $0<\epsilon\leq 3$. In this case,  by
\cite[pp. 73-76]{DGF}, we know that, for $d>30$, $S$ is contained in
a smooth rational normal scroll  of $\mathbb P^5$ of dimension $3$.
Taking into account that we may assume $K^2_S\leq 8\chi (\mathcal
O_S)$ (compare with  Remark \ref{k8}, $(i)$ and $(ii)$), at this
point Theorem \ref{lbound} follows from \cite[Proposition 2.2]{DGF},
when $d>30$.

\smallskip
In order to examine the remaining cases $26\leq d \leq 30$, we
refine the analysis appearing in \cite{DGF}. In fact, we are going
to prove that, assuming $r=5$, $S$ is a scroll, and (\ref{bound}),
it follows that $S$ is contained in a smooth rational normal scroll
of $\mathbb P^5$ of dimension $3$ also when $26\leq d \leq 30$. Then
we may conclude as before, because \cite[Proposition 2.2]{DGF} holds
true for $d\geq 18$.

\smallskip
First, observe that if $S$ is contained in a threefold $T\subset
\mathbb P^5$ of dimension $3$ and minimal degree $3$, then $T$ is
necessarily a {\it smooth} rational normal scroll \cite[p. 76]{DGF}.
Moreover, observe that we may apply the same argument as in \cite[p.
75-76]{DGF} in order to exclude the case  $S$ is contained in a
threefold of degree $4$. In fact the argument works for $d>24$
\cite[p. 76, first line after formula (13)]{DGF}.

\smallskip
In conclusion, assuming $r=5$, $S$ is a scroll, and (\ref{bound}),
it remains to exclude that $S$ is not contained in a threefold of
degree $<5$, when $26\leq d \leq 30$.

\smallskip
Assume $S$ is not contained in a threefold of degree $<5$. Denote by
$\Gamma\subset \mathbb P^3$ the general hyperplane section of $H$.
Recall that $26\leq d \leq 30$.

\bigskip $\bullet$ Case I: $h^0(\mathbb P^3,\mathcal
I_{\Gamma}(2))\geq 2$.

It is impossible. In fact, if $d>4$, by monodromy \cite[Proposition
2.1]{CCD}, $\Gamma$ should be contained in a reduced and irreducible
space curve of degree $\leq 4$, and so, for $d>20$, $S$ should be
contained in a threefold of degree $\leq 4$  \cite[Theorem
(0.2)]{CC}.

\bigskip $\bullet$ Case II: $h^0(\mathbb P^3,\mathcal
I_{\Gamma}(2))=1$ and $h^0(\mathbb P^3,\mathcal I_{\Gamma}(3))>4$.

As before, if $d>6$, by monodromy, $\Gamma$ is contained in a
reduced and irreducible space curve  $X$ of degree $\deg(X)\leq 6$.
Again as before, if $\deg(X)\leq 4$, then $S$ is contained in a
threefold of degree $\leq 4$. So we may assume $5\leq \deg(X)\leq
6$.

Since $d\geq 26$, by Bezout's Theorem we have $h_{\Gamma}(i)=h_X(i)$
for all $i\leq 4$. Let $X'$ be the general plane section of $X$.
Since $h_X(i)\geq \sum_{j=0}^{i}h_{X'}(j)$, we have $h_X(3)\geq 14$
and $h_X(4)\geq 19$ \cite[pp. 81-87]{EH}. Therefore, when $d\geq
26$, taking into account \cite[Corollary (3.5)]{EH}, we get:
$$
h_{\Gamma}(1)=4,\, h_{\Gamma}(2)=9,\, h_{\Gamma}(3)\geq 14,\,
h_{\Gamma}(4)\geq 19,
$$
$$
h_{\Gamma}(5)\geq 22, \,
h_{\Gamma}(6)\geq \min\{d,\, 27\},\, h_{\Gamma}(7)=d.
$$

It follows that:
$$
p_a(C)\leq \sum_{i=1}^{+\infty}d-h_{\Gamma}(i)\leq
(d-4)+(d-9)+(d-14)+(d-19)+(d-22)+3=5d-65,
$$
which is $<\frac{1}{8}d(d-6)+1$ for $d \geq 26$. This is in contrast
with (\ref{bound}).

\bigskip $\bullet$ Case III: $h^0(\mathbb P^3,\mathcal
I_{\Gamma}(2))=1$ and $h^0(\mathbb P^3,\mathcal I_{\Gamma}(3))=4$.

We have:
$$
h_{\Gamma}(1)=4,\, h_{\Gamma}(2)=9,\, h_{\Gamma}(3)=16, \,
h_{\Gamma}(4)\geq 19,
h_{\Gamma}(5)\geq 24, \, h_{\Gamma}(6)=d.
$$

It follows that:
$$
p_a(C)\leq \sum_{i=1}^{+\infty}d-h_{\Gamma}(i)\leq
(d-4)+(d-9)+(d-16)+(d-19)+(d-24)=5d-72,
$$
which is $< \frac{1}{8}d(d-6)+1$ for $d \geq 26$. This is in
contrast with (\ref{bound}).

\bigskip $\bullet$ Case IV: $h^0(\mathbb P^3,\mathcal
I_{\Gamma}(2))=0$.

We have:
$$
h_{\Gamma}(1)=4,\, h_{\Gamma}(2)=10,\, h_{\Gamma}(3)\geq 13, \,
h_{\Gamma}(4)\geq 19,
$$
$$
h_{\Gamma}(5)\geq 22,\, h_{\Gamma}(6)\geq \min\{d,\, 28\}, \,
h_{\Gamma}(7)=d.
$$
It follows that:
$$
p_a(C)\leq  \sum_{i=1}^{+\infty}d-h_{\Gamma}(i)\leq
(d-4)+(d-10)+(d-13)+(d-19)+(d-22)+2=5d-66,
$$
which is $< \frac{1}{8}d(d-6)+1$ for $d \geq 26$. This is in
contrast with (\ref{bound}).

This concludes the proof of Theorem \ref{lbound}.

\begin{remark}\label{altro}
$(i)$ Let $Q\subseteq \mathbb P^3$ be a smooth quadric, and
$H\in|\mathcal O_Q(1,d-1)|$ be a smooth rational curve of degree $d$
\cite[p. 231, Exercise 5.6]{Hartshorne}. Let $S\subseteq\mathbb P^4$
be the projective cone over $H$. A computation, which we omit,
proves that
$$
\chi (\mathcal O_S)=1-\binom{d-1}{3}.
$$
Therefore, if $S$ is singular, it may happen that $\chi (\mathcal
O_S)<-\frac{1}{8}d(d-6)$. One may ask whether $1-\binom{d-1}{3}$ is
a lower bound for $\chi(\mathcal O_S)$ for every {\it integral}
surface.

\smallskip
$(ii)$ Let $(S,\mathcal L)$ be a smooth surface, polarized by a very
ample line bundle $\mathcal L$ of degree $d$. By Harris' bound for
the geometric genus $p_g(S)$ of $S$ \cite{H}, we see that
$p_g(S)\leq \binom{d-1}{3}$. Taking into account that for a smooth
surface one has $\chi(\mathcal O_S)=h^0(S,\mathcal
O_S)-h^1(S,\mathcal O_S)+h^2(S,\mathcal O_S) \leq 1+h^2(S,\mathcal
O_S)=1+p_g(S)$, from Theorem \ref{lbound} we deduce (the first
inequality only when $d>25$):
$$
-\binom{\frac{d}{2}-1}{2}\leq \chi (\mathcal O_S)\leq
1+\binom{d-1}{3}.
$$
\end{remark}

\bigskip
\section{Proof of Corollary \ref{Veronese}}

$\bullet$ First, assume that $\chi(\mathcal
O_S)=-\frac{1}{8}d(d-6)$.

\smallskip
By Theorem \ref{lbound}, we know that $r=5$. Moreover, $S$ is
contained in a nonsingular threefold $T\subseteq \mathbb P^5$ of
minimal degree $3$. Therefore, the general hyperplane section
$H=S\cap L$ of $S$ ($L\cong \mathbb P^4$ denotes the general
hyperplane of $\mathbb P^5$) is contained in a smooth surface
$\Sigma=T\cap L$ of $L\cong \mathbb P^4$, of minimal degree $3$.

This surface $\Sigma$ is isomorphic to the blowing-up of $\mathbb
P^2$ at a point, and, for a suitable point $x\in V\backslash L$, the
projection of $\mathbb P^5\backslash\{x\}$ on $L\cong \mathbb P^4$
from $x$ restricts to an isomorphism
$$p_x:V\backslash\{x\}\to \Sigma\backslash E,$$ where $E$ denotes
the exceptional line of $\Sigma$ \cite[p. 58]{BV}.

Since $S$ is linearly equivalent on $T$ to $\frac{d}{2}(H_T- W_T)$
($H_T$ denotes the hyperplane section of $T$, and $W_T$ the ruling),
it follows that $H$ is linearly equivalent on $\Sigma$ to
$\frac{d}{2}(H_{\Sigma}- W_{\Sigma})$ (now $H_{\Sigma}$ denotes the
hyperplane section of $\Sigma$, and $W_{\Sigma}$ the ruling of
$\Sigma$). Therefore, $H$ does not meet the exceptional line
$E=H_{\Sigma}- 2W_{\Sigma}$. In fact, since $H_{\Sigma}^2=3$,
$H_{\Sigma}\cdot W_{\Sigma}=1$, and $W_{\Sigma}^2=0$, one has:
$$(H_{\Sigma}- W_{\Sigma})\cdot (H_{\Sigma}-
2W_{\Sigma})= H_{\Sigma}^2-3H_{\Sigma}\cdot
W_{\Sigma}+2W_{\Sigma}^2=0.$$

This implies that $H$ is contained in $\Sigma\backslash E$, and the
assertion of Corollary \ref{Veronese} follows.

\smallskip
$\bullet$ Conversely, assume there exists a curve $C$ on the
Veronese surface $V\subseteq \mathbb P^5$, and a point $x\in
V\backslash C$, such that $H$ is the projection $p_x(C)$ of $C$ from
the point $x$.

\smallskip
In particular, $d$ is an even number, and $H$ is contained in a
smooth surface $\Sigma\subseteq L\cong \mathbb P^4$ of minimal
degree, and is disjoint from the exceptional line $E\subseteq
\Sigma$. By \cite[Theorem (0.2)]{CC}, $S$ is contained in a
threefold $T\subseteq \mathbb P^5$ of minimal degree. $T$ is
nonsingular. In fact, otherwise, $H$ should be a Castelnuovo's curve
in $\mathbb P^4$ \cite[p. 76]{DGF}. On the other hand, by our
assumption, $H$ is isomorphic to a plane curve of degree
$\frac{d}{2}$. Hence, we should have:
$$
g=\frac{d^2}{6}-\frac{2}{3}d+1=\frac{d^2}{8}-\frac{3}{4}d+1
$$
(the first equality because $H$ is Castelnuovo's, the latter because
$H$ is isomorphic to a plane curve of degree $\frac{d}{2}$). This is
impossible when $d>0$.

Therefore, $S$ is contained in a smooth threefold $T$ of minimal
degree in $\mathbb P^5$.

Now observe that in $\Sigma$ there are only two families of curves
of degree even $d$ and genus $g=\frac{d^2}{8}-\frac{3}{4}d+1$. These
are the curves linearly equivalent on $\Sigma$ to
$\frac{d}{2}(H_{\Sigma}- W_{\Sigma})$, and the curves equivalent to
$\frac{d+2}{6}H_{\Sigma}+ \frac{d-2}{2}W_{\Sigma}$. But only in the
first family the curves do not meet $E$. Hence, $H$ is linearly
equivalent on $\Sigma$ to $\frac{d}{2}(H_{\Sigma}- W_{\Sigma})$.
Since the restriction ${\text{Pic}}(T)\to {\text{Pic}}(\Sigma)$ is
bijective, it follows that $S$ is linearly equivalent on $T$ to
$\frac{d}{2}(H_{T}- W_{T})$. By Theorem \ref{lbound}, $S$ is  a
fortiori linearly normal, and of minimal Euler characteristic
$\chi(\mathcal O_S)=-\frac{1}{8}d(d-6)$.

\end{document}